\documentclass[12pt]{iopart}
\usepackage{iopams}

\newtheorem{thm}{Theorem}[section]

\begin{document}

\title[]{Contact Equivalence Problem for KDV-type equations}

\author{Mostafa Hesamiarshad and Mehdi Nadjafikhah}

\address{Department of Mathematics, College of Basic Sciences, Karaj Branch, Islamic Azad University,
              Alborz, Iran}
\ead{m.hesami@tuyiau.ac.ir and m$_{-}$nadjafikhah@iust.ac.ir}
\begin{abstract}
The Cartan's method of equivalence and moving coframe method has
been applied to solve the local equivalence problem for KDV-type
equations under the action of a pseudo-group of contact
transformations. The structure equations, the sets of differential
invariants for symmetry groups and equivalent
conditions of these equations are found.
\end{abstract}

\maketitle
\section*{Introduction}
In the beginning of twentieth century, Elie Cartan developed a
uniform method for analyzing the differential invariants of many
geometric structures, nowadays called the `Cartan equivalence
method'. Also, the method of equivalence is a systematic procedure
that allows one to decide whether two systems of differential
equations can be mapped one to another by a transformation taken
in a given pseudo-group. Later, C. Erhesmann and S. Chern
introduced two important concepts to the method of equivalence:
jets spaces and G-structures. In recent years, with the help of
mathematical software, many authors have successfully applied the
method of equivalence to many interesting problems:
classifications of differential equations \cite{FO1,KNM,KNT},
holonomy groups \cite{RBT}, inverse variational problems
\cite{FME} and general relativity \cite{BCI,KAA}.

In this paper, we consider a local equivalence problem for the
class of equations
   \begin{equation}\label{a}
u_{xxx}=u_t+Q(u,u_x)
\end{equation}
under the contact transformation of a pseudo-group. Equation
(\ref{a}) is the standard Korteweg-de Vries (KDV)
 equation if $Q=uu_x$, it is the Modified KDV-equation if $Q=u^2u_x$, and
 it is generalized KDV equation if $Q=h(u)u_x$.
Two equations are said to be equivalent if there exists a contact
transformation maps one equation to another. We use Elie Cartan's
method of equivalence \cite{EC,EIS,TME} which in form developed by
Fels and Olver \cite{FO1,FO2} and as stated by morozov \cite{MC}
to compute the Maurer - Cartan forms, structure equations and
basic invariants for symmetry groups of equations.
 Cartan's solution to the
equivalence problem states that two submanifolds are (locally)
equivalent if and only if their classifying manifolds (locally)
overlap. The symmetry classification problem for classes of
differential equations is closely related to the problem of local
equivalence: symmetry groups of two equations are necessarily
isomorphic if these equations are equivalent, while in general the
converse of this issue is not true. For the symmetry analysis of
(\ref{a}) the reader is referred to \cite{SK}.

\section{Pseudo-group of contact transformations of differential
equations } In this section we describe the local equivalence
problem for differential equations under the action of the pseudo
group of contact transformations. Two equations are said to be
equivalent if there exists a contact transformation which maps the
equations to each other. We apply Elie Cartan's structure theory
of Lie pseudo-groups  to obtain necessary and sufficient
conditions under which equivalence mappings can be found. This
theory describes a Lie pseudo-group in terms of a set of invariant
differential 1-forms called Maurer-Cartan forms, which contain all
information about the seudo-group. In particular, they give basic
invariants and operators of invariant differential, which in terms
allow us to solve equivalence problem for submanifolds under the
action of the pseudo-group. Recall that expansions of exterior
differentials of Maurer-Cartan forms in terms of the form
themselves, yields the Cartan structure equation for the
prescribed pseudo-group.
\\
Suppose
$\pi:\mathbb{R}^n\times\mathbb{R}^m\rightarrow\mathbb{R}^n$ is a
trivial bundle with the local base coordinates $(x_1,\ldots,x_n)$
and the local fibre coordinates $(u_1,\ldots,u_m)$; then
$J^1(\pi)$ is denoted by the bundle of the first-order jets of
sections of $\pi$, with the local coordinates $(x_i, u_{\alpha},
p^{\alpha}_i ), i\in\{1, \ldots, n\}, {\alpha}\in \{1,\ldots,m\}$,
where $p^{\alpha}_i=\frac{\partial u_{\alpha}}{\partial x_i}$. For
every local section $\left(x_i, f_{\alpha}\left(x\right)\right)$
of $\pi$, the corresponding 1-jet
$\left(x_i,f_{\alpha}\left(x\right),\frac{\partial
f_{\alpha}\left(x\right)}{\partial x_i}\right)$ is denoted by
$j_1(f)$. A differential 1-form $\nu$ on $J^1(\pi)$ is called a
contact form, if it is annihilated by all 1-jets of local
sections: $j_1(f)^*\nu=0$. In the local coordinates every contact
1-form is a linear combination of the forms
$\nu^{\alpha}=du_{\alpha}-p^{\alpha}_idx_i,
{\alpha}\in\{1,\ldots,m\}$. A local diffeomorphism:
\begin{equation}\label{delta}
\Delta:J^1(\pi)\rightarrow J^1(\pi),\qquad
\Delta:(x_i,u_{\alpha},p^{\alpha}_i)\rightarrow
(\overline{x}_i,\overline{u}_{\alpha},\overline{p}^{\alpha}_i)
\end{equation}
is called a contact transformation, if for every contact 1-form
$\nu$, the form $\Delta^*\overline{\nu}$ is also contact. Suppose
$\mathcal{R}$ is a first-order differential equation in $m$
dependent and $n$ independent variables. We consider $\mathcal{R}$
as a sub-bundle in $J^1(\pi)$. Suppose $Cont(\mathcal{R})$ is the
group of contact symmetries for $\mathcal{R}$. It consists of all
the contact transformations on $J^1(\pi)$ mapping $\mathcal{R}$ to
itself.

It was shown in \cite{MC} that the following differential 1-forms:
\begin{eqnarray}
  \Theta^{\alpha} &=& a^{\alpha}_{\beta}(du_{\beta}-p^{\beta}_jdx_j), \nonumber\\
  \Xi^i &=& b^i_jdx_j+c^i_{\beta}\Theta^{\beta}, \nonumber\\
  \Sigma^{\alpha}_i &=& a^{\alpha}_{\beta}B^i_jdp^{\beta}_j+f^{\alpha}_{i\beta}\Theta^{\beta}+g^{\alpha}_{ij}\Xi^j.\nonumber
    \end{eqnarray}
 are the Maurer-Cartan forms of $Cont(J^1(\pi))$. They are defined on $J^1(\pi) \times \mathcal{H}$, where
$\mathcal{H}={(a^{\alpha}_{\beta},b^{i}_{j},c^{i}_{\beta},f^{\alpha}_{i\beta},g^{\alpha}_{ij})}
\mid \alpha,\beta \in \{1,\ldots,m\},i,j\in \{1,\ldots,n\},\\
det(a^{\alpha}_{\beta}).det(b^i_j)\neq0,
g^{\alpha}_{ij}=g^{\alpha}_{ji}$ and $(B^i_j)$ is the inverse
matrix for $(b^i_j)$. They satisfy the structure equations
\begin{eqnarray}
  d\Theta^{\alpha} &=& \Phi^{\alpha}_{\beta} \wedge \Theta^{\beta}+\Xi^k \wedge \Sigma^{\alpha}_k, \nonumber\\
  d\Xi^i &=&\Psi^{i}_{k} \wedge \Xi^{k}+\Pi^{i}_{\gamma} \wedge \Theta^{\gamma}, \nonumber\\
  d\Sigma^{\alpha}_i &=& \Phi^{\alpha}_{\gamma} \wedge
  \Sigma^{\gamma}_i-\Psi^k_i \wedge
  \Sigma^{\alpha}_k+\Lambda^{\alpha}_{i\beta} \wedge
  \Theta^{\beta}+\Omega^{\alpha}_{ij} \wedge \Xi^j,\nonumber
\end{eqnarray}
where the forms
$\Phi^{\alpha}_{\beta},\Psi^{i}_{j},\Pi^{i}_{\beta},\Lambda^{\alpha}_{i\beta}$
and $\Omega^{\alpha}_{ij}$ depend on differentials of the
coordinates of $\mathcal{H}$.\\

Differential equations defines a submanifold $\mathcal{R}\subset
J^1(\pi)$. The Maurer-Cartan forms for its symmetry pseudo-group
$Cont(\mathcal{R})$ can be found from restrictions
$\theta^{\alpha}=\imath^*\Theta^{\alpha}, \xi^i=\imath^*\Xi^i$and
$\sigma^{\alpha}_i=\imath^*\Sigma^{\alpha}_i$, where
$\imath=\imath_0 \times id
:\mathcal{R}\times\mathcal{H}\longrightarrow
 J^1(\pi)\times\mathcal{H}$ with
$\imath_0:\mathcal{R} \longrightarrow J^1(\pi)$ is defined by our
differential equations. In order to compute the Maurer-Cartan
forms for the symmetry pseudo-group, we implement Cartan's
equivalence method. Firstly, the forms $\theta^{\alpha},
\xi^i,\sigma^{\alpha}_i$ are linearly dependent, i.e. there exists
a nontrivial set of functions $U_{\alpha}, V_i, W^i_{\alpha}$ on
$\mathcal{R}\times\mathcal{H}$ such that
$U_{\alpha}\theta^{\alpha}+ V_i\xi^i+
W^i_{\alpha}\sigma^{\alpha}_i\equiv0$. Setting these functions
equal to some appropriate constants allows us to introduce a part
of the coordinates of $\mathcal{H}$ as functions of the other
coordinates of $\mathcal{R}\times\mathcal{H}$. Secondly, we
substitute the obtained values into the forms
$\phi^{\alpha}_{\beta}=\imath^*\Phi^{\alpha}_{\beta}$ and
$\psi^i_k=\imath^*\Psi^i_k$ coefficients of semi-basic forms
$\phi^{\alpha}_{\beta}$ at $\sigma^{\gamma}_j, \xi^j$, and the
coefficients of semi-basic forms $\psi^{i}_{j}$ at
$\sigma^{\gamma}_j$ are lifted invariants of $Cont(\mathcal{R})$.
We set them equal to appropriate constants and get expressions for
the next part of the coordinates of $\mathcal{H}$, as functions of
the other coordinates of $\mathcal{R}\times\mathcal{H}$. Thirdly,
we analyze the reduced structure equations
\begin{eqnarray}
  d\theta^{\alpha} &=& \phi^{\alpha}_{\beta} \wedge \theta^{\beta}+\xi^k \wedge \sigma^{\alpha}_k, \nonumber\\
  d\xi^i &=&\psi^{i}_{k} \wedge \xi^{k}+\pi^{i}_{\gamma} \wedge \theta^{\gamma}, \nonumber\\
  d\sigma^{\alpha}_i &=& \phi^{\alpha}_{\gamma} \wedge
  \sigma^{\gamma}_i-\psi^k_i \wedge
  \sigma^{\alpha}_k+\lambda^{\alpha}_{i\beta} \wedge
  \theta^{\beta}+\omega^{\alpha}_{ij} \wedge \xi^j.\nonumber
\end{eqnarray}
If the essential torsion coefficients are dependent on the group
parameters , then we may normalize them to constants and find some
new part of the group parameters, which, upon being substituted
into the reduced modified Maurer-Cartan forms, allows us to repeat
the procedure of normalization. This process has two results.
First, when the reduced lifted coframe appears to be involutive,
this coframe is the desired set of defining forms for
$Cont(\mathcal{R})$. Second, when the coframe is not involutive we
should apply the procedure of prolongation described in
\cite{EIS}.
\section{Structure and invariants of symmetry groups for  KDV-type equations }
 Consider the following
 system equivalent to (\ref{a}) of first order:
\begin{equation}\label{P}
    u_x=v \qquad  v_{x}=w,\qquad  w_{x}=u_t+Q(u,v).
\end{equation}
We apply the method described in the previous section to the class
of equations (\ref{P}).
 We denotes that  $ t=x_1,x=x_2, u=u_1, v=u_2, w=u_3, u_t=p^1_1,
u_x=p^1_2, v_t=p^2_1, v_x=p^2_2, w_t=p^3_1, w_x=p^3_2$. We
consider this system as a sub-bundle of the bundle
$J^1(\pi),\pi:\mathbb{R}^2\times\mathbb{R}^3\longrightarrow\mathbb{R}^2$,
with local coordinates $\{x_1, x_2, u_1, u_2,u_3, p^1_1, p^2_1,
p^3_1\}$, where the embedding $\iota$ is defined by the following
equalities:
\begin{equation}\label{EQ}
  p^1_2=u_2 \qquad p^2_2=u_3 \qquad
  p^3_2=p^1_1+Q(u_1,u_2)
\end{equation}
The forms $\theta^{\alpha}=\iota^*\Theta^{\alpha},\alpha \in
\{1,2,3\},\xi^{i}=\iota^*\Xi^{i},i \in \{1,2\}$, are linearly
dependent. The group parameters $a^{\alpha}_{\beta},b^i_j$ should
satisfy the simultaneous  conditions
$det(a^{\alpha}_{\beta})\neq0,det(b^i_j)\neq0$.
 Linear dependence between the forms $\sigma^{\alpha}_i$ are
\begin{equation}\label{DEP}
 \sigma^1_2=0 , \qquad \sigma^2_2=0, \qquad \sigma^3_2=\sigma^1_1
\end{equation}
Computing the linear dependence conditions (\ref{DEP}) gives the
group parameters $a^1_1, a^1_2, a^1_3,\\ a^2_3, b^1_2, f^1_{11},
f^1_{21}, f^1_{22},f^1_{23}, f^2_{21},f^2_{22}, f^2_{23},f^3_{22},
f^3_{23}, g^1_{12}, g^1_{22}, g^2_{12}, g^2_{22}, g^3_{12},
g^3_{22}$ as functions of other group parameters and the local
coordinates $\{x_1,x_2,u_1,u_2,u_3æ p^1_1,p^2_1,p^3_1\}$ of
$\mathcal{R}$.In particular,
\begin{eqnarray}
&&f^2_{21}=-\frac{a^2_1a^3_2a^2_2-{(a^2_1)}^2a^3_2-a^3_1{(a^2_2)}^2+g^2_{22}c^2_1{(a^3_3)}^2b^1_1a^2_2+g^2_{12}c^1_1{(a^3_3)}^2b^1_1a^2_2}{{(a^3_3)}^2b^1_1a^2_2},\nonumber\\
&&g^2_{12}=\frac{-p^3_1a^2_2b^2_2-p^2_1a^2_1b^2_2+Qb^2_1a^2_2+p^1_1b^2_1a^2_2+u_3b^2_1a^2_1+u_3b^2_1a^2_2}{b^1_1{(b^2_2)}^2}, \hspace{.8cm} b^1_2=0,\nonumber\\
&&f^2_{22}=-\frac{a^3_2a^2_2+a^2_1a^3_3+g^2_{22}c^2_2a^2_2b^2_2a^3_3+g^2_{12}c^1_2a^2_2b^2_2a^3_3}{a^3_3b^2_2a^2_2}, \hspace{2.9cm} a^2_3=0,\nonumber\\
&&f^1_{22}=-\frac{a^3_3b^1_1+g^1_{22}c^2_2a^2_2{(b^2_2)}^2+g^1_{12}c^1_2a^2_2{(b^2_2)}^2}{a^2_2{(b^2_2)}^2}, \hspace{3.9cm} a^1_3=0,\nonumber\\
&&f^1_{21}=\frac{a^2_1-g^1_{22}c^2_1a^2_2b^2_2-g^1_{12}c^1_1a^2-2b^2_2}{a^2_2b^2_2}, \hspace{4.8cm} a^1_2=0, \nonumber\\
&&f^2_{23}=-\frac{a^2_2+g^2_{22}c^2_3a^3_3b^2_2+g^2_{12}c^1_3a^3_3b^2_2}{a^3_3b^2_2},\hspace{4.6cm}a^1_1=\frac{a^3_3b^1_1}{b^2_2} ,\nonumber\\
&&g^2_{22}=-\frac{Qa^2_2+p^1_1a^2_2+u_3a^2_1+u_3a^2_2}{{(b^2_2)}^2}, \hspace{2.7cm} f^1_{23}=-g^1_{22}c^2_3-g^1_{12}c^1_3,\nonumber\\
&&g^1_{12}=\frac{a^3_3(u_3b^2_1-p^2_1b^2_2)}{{(b^2_2)}^2},
\hspace{5.2cm}g^1_{22}=\frac{u_3a^3_3b^1_1}{{(b^2_2)}^2}.\nonumber
\end{eqnarray}

The expressions for $f^1_{11}, f^3_{22}, f^3_{23} ,g^3_{12}$ and
$g^3_{22}$ are too long to be written out here completely.\\ The
analysis of the semi-basic modified Maurer-Cartan forms
$\phi^{\alpha}_{\beta}, \psi^i_k$ at the obtained values of the
group parameters gives the following
normalization.\\
 The form $\psi^1_2$ is semi-basic, and
$\psi^1_2\equiv-c^1_3\sigma^1_1$. So we take $c^1_3=0$. For the
semi-basic form $\phi^1_2$ we have
\begin{eqnarray}
    \phi^1_2&\equiv&c^1_2\sigma^1_1 \quad(mod  \quad \theta^1,\theta^2,\theta^3,\xi^1,\xi^2),\qquad\qquad\qquad\qquad\qquad \nonumber
  \end{eqnarray}
  thus, we can assume $c^1_2=0$. And for the
semi-basic form $\phi^1_3$ we have
  \begin{eqnarray}
  \phi^1_3&\equiv&-\frac{c^2_3a^3_3b^2_1u_3-c^2_3a^3_3b^2_2p^2_1+f^1_{13}(b^2_2)^3}{(b^2_2)^3}\xi^1 \quad(mod  \quad \theta^1,\theta^2,\theta^3),\nonumber
\end{eqnarray}
so we set the coefficient at $\xi^1$ equal to 0 and find
\begin{eqnarray}
f^1_{13}=\frac{c^2_3a^3_3(b^2_2p^2_1-b^2_1u_3)}{(b^2_2)^3}.\qquad\qquad\qquad\qquad\qquad\qquad\qquad\qquad\nonumber
\end{eqnarray}

By doing the analysis of the modified semi-basic Maurer-Cartan
forms in the same way, we can normalize the following group
parameters:
\begin{eqnarray}
&&a^2_1=a^3_1=a^3_2=0, \hspace{6cm} a^2_2=a^3_3b^2_2,\nonumber \\
&&b^1_1={(b^2_2)}^3, \hspace{7cm} b^2_1=-b^2_2Q_v,\nonumber\\
&&c^1_1=c^2_1=c^2_2=c^2_3=0,\nonumber\\
&&f^1_{12}=f^2_{13}=f^3_{21}=f^3_{12}=f^2_{12}=f^3_{11}=0,\hspace{1.1cm}f^3_{13}=f^2_{12}=\frac{Q_u}{{(b^2_2)}^3}.\nonumber
\end{eqnarray}

 We denote by $\mathcal{S}_1$ the
subclass of equations (\ref{P}) such that: $Q_{u^2}= 0$, $Q_{uv}=
0$ and $Q_{v^2}= 0$. In this case, the structure equations of the
symmetry group for system (\ref{P}) have the form:
\begin{eqnarray}\label{ES2}
   && d\theta^1=-\theta^1 \wedge( \eta_4+2\eta_5)-\theta^2 \wedge \xi^2+\xi^1 \wedge\sigma^1_1,\nonumber\\
   && d\theta^2=-\theta^2 \wedge (\eta_4+\eta_5)-\theta^3\wedge \xi^2+\xi^1 \wedge\sigma^2_1,\nonumber\\
   && d\theta^3=-\theta^3 \wedge \eta_4+\xi^1 \wedge\sigma^3_1+\xi^2 \wedge\sigma^1_1,\nonumber\\
    && d\xi^1=-3\xi^1 \wedge \eta_5,\\
    && d\xi^2=-\xi^2 \wedge \eta_5,\nonumber\\
    && d\sigma^1_1=-\xi^1 \wedge \eta_3+\xi^2 \wedge\sigma^2_1-\sigma^1_1\wedge (\eta_4-\eta_5),\nonumber\\
     && d\sigma^2_1=-\xi^1 \wedge \eta_1+\xi^2 \wedge\sigma^3_1-\sigma^2_1\wedge (\eta_4-2\eta_5),\nonumber\\
    && d\sigma^3_1=-\xi^1 \wedge \eta_2-\xi^2 \wedge\eta_3+\sigma^3_1\wedge (\eta_4-3\eta_5).\nonumber
   \end{eqnarray}

In the analysis of structure equations we can't absorb any group
parameters more than before. Besides, the Cartan character is
$s_1=5$ and the indetermination degree is 3, thus the involution
test fails. So we adopt the procedures of prolongation to compute
the new structure equations:
\begin{eqnarray}\label{ES1}
   && d\theta^1=-\theta^1 \wedge( \eta_4+2\eta_5)-\theta^2 \wedge \xi^2+\xi^1 \wedge\sigma^1_1,\nonumber\\
   && d\theta^2=-\theta^2 \wedge (\eta_4+\eta_5)-\theta^3\wedge \xi^2+\xi^1 \wedge\sigma^2_1,\nonumber\\
   && d\theta^3=-\theta^3 \wedge \eta_4+\xi^1 \wedge\sigma^3_1+\xi^2 \wedge\sigma^1_1,\nonumber\\
    && d\xi^1=-3\xi^1 \wedge \eta_5,\nonumber\\
    && d\xi^2=-\xi^2 \wedge \eta_5,\\
    && d\sigma^1_1=-\xi^1 \wedge \eta_3+\xi^2 \wedge\sigma^2_1-\sigma^1_1\wedge (\eta_4-\eta_5),\nonumber\\
     && d\sigma^2_1=-\xi^1 \wedge \eta_1+\xi^2 \wedge\sigma^3_1-\sigma^2_1\wedge (\eta_4-2\eta_5),\nonumber\\
    && d\sigma^3_1=-\xi^1 \wedge \eta_2-\xi^2 \wedge\eta_3-\sigma^3_1\wedge (\eta_4-3\eta_5),\nonumber\\
     && d\eta_1=-\beta_1\wedge\xi^1+\xi^2 \wedge \eta_2-\eta_1\wedge(\eta_4-5\eta_5),\nonumber\\
    && d\eta_2=-\beta_2\wedge\xi^1-\beta_3 \wedge \xi^2-\eta_2\wedge(\eta_4-6\eta_5),\nonumber\\
    && d\eta_3=-\beta_3\wedge\xi^1+\xi^2 \wedge \eta_1-\eta_3\wedge(\eta_4-4\eta_5),\nonumber\\
    && d\eta_4=0,\nonumber\\
    && d\eta_5=0.\nonumber
\end{eqnarray}
In structure equations (\ref{ES1}), the forms
$\eta_1,\cdots\eta_5$ on $J^2(\pi) \times \mathcal{H}$ depend on
differentials of the parameters of $\mathcal{H}$, while the forms
$\beta_1,\beta_2,\beta_3$ depend on differentials of the
prolongation variables. In the structure equations (\ref{ES1}) the
degree of indetermination is 3 and the Cartan characters are
$s_1=3,s_2= \ldots = s_{13} = 0$. Consequently,  Cartan's test for
the lifted coframe
$\{\theta^1,\theta^2,\theta^3,\xi^1,\xi^2,\sigma^1_1,\sigma^2_1,\sigma^3_1,\eta_1,\eta_2,\eta_3,\eta_4,\eta_5\}$
is satisfied. Therefore, the coframe is involutive. All the
essential torsion coefficients in the structure equations
(\ref{ES1}) are constant. By applying Theorem 11.8 of \cite{EIS},
we have:
\begin{thm}
All systems from $\mathcal{S}_1$ are (locally) equivalent under
contact transformations.
\end{thm}

We denote by $\mathcal{S}_2$ the subclass of equations (\ref{P})
such that $Q_{u^2}=0, Q_{v^2}= 0$ and $Q_{uv}\neq 0$.
$\mathcal{S}_2$  is is completely described by equations of the
form $u_{xxx}=u_t+Au+Bu_x+Cuu_x+D$ with $A,B,C,D \in \mathbb{R}$
and $C \neq 0$. KDV equation belongs to $\mathcal{S}_2$. The
analysis of the structure equations gives the following essential
torsion coefficients and the corresponding normalization:
\begin{eqnarray}
d\sigma^2_1=
\omega^2_{11}\wedge\xi^1+\phi^3_3\wedge\sigma^2_1-2\psi^2_2
\wedge\sigma^2_1+\frac{C}{{(b^2_2)}^4a^3_3}\theta^1\wedge\theta^3+\cdots\nonumber
\end{eqnarray}
thus, we can assume $a^3_3=\frac{C}{{(b^2_2)}^4}$; and
 \begin{eqnarray}
d\sigma^2_1=\omega^3_{11}\wedge\xi^2-6\omega^1_{11}\wedge\sigma^1_1+\theta^1\wedge\theta^2-\frac{Cv}{(b^2_2)^3}\theta^3\wedge\xi^2+\cdots\nonumber
\end{eqnarray}
thus, we can assume $b^2_2=\sqrt[3]{Cv}$ .

After this normalization, the structure equations of coframe
$\{\theta^1,\theta^2,\theta^3,\xi^1,\xi^2,\\\sigma^1_1,\sigma^2_1,\sigma^3_1\}$
is,
\begin{eqnarray}\label{ES2}
&&d\theta^1=\frac{2}{3}\theta^1\wedge\theta^2+(\frac{2}{3}I_2-I_3-1)\theta^1\wedge\xi^1+\frac{2}{3}I_1\theta^1\wedge\xi^2-\theta^2\wedge\xi^2+\xi^1\wedge\sigma^1_1,\nonumber\\
&&d\theta^2=(I_3-I_2-1)theta^2\wedge\xi^1+I_1\theta^2\wedge\xi^2-\theta^3\wedge\xi^2+\xi^1\wedge\sigma^2_1 ,\nonumber\\
&&d\theta^3=-\frac{4}{3}\theta^2\wedge\theta^3-(\frac{4}{3}I_2-I_3-1)\theta^3\wedge\xi^1+\frac{4}{3}I_1\theta^3\wedge\xi^2+\xi^2\wedge\sigma^1_1+\xi^1\wedge\sigma^3_1 ,\nonumber\\
&&d\xi^1=\theta^2\wedge\xi^1-I_1\xi^1\wedge\xi^2 ,\\
&&d\xi^2=-\theta^1\wedge\xi^1+\frac{1}{3}\theta^2\wedge\xi^21(\frac{1}{2}I_2-1)\xi^1\wedge\xi^2 ,\nonumber\\
&&d\sigma^1_1=\pi_3\wedge\xi^1-I_1\theta^2\wedge\xi^2-\theta^3\wedge\xi^2-2\theta^2\wedge\sigma^2_1-2I_1\xi^2\wedge\sigma^2_1\nonumber\\
&& \qquad \quad+\xi^2\wedge\sigma^3_1 ,\nonumber\\
&&d\sigma^2_1=\pi_1\wedge\xi^1+\theta^1\wedge\theta^3-I_1\theta^2\wedge\xi^2-\theta^3\wedge\xi^2-2\theta^2\wedge\sigma^2_1
-2I_1\xi^2\wedge\sigma^2_1\nonumber\\
&& \qquad \quad+\xi^2\wedge\sigma^3_1 ,\nonumber\\
&&d\sigma^3_1=\pi_2\wedge\xi^1+\pi_3\wedge\xi^2+\theta^1\wedge\sigma^1_1+\theta^2\wedge\theta^3-I_1\theta^3\wedge\xi^2-\frac{7}{3}\theta^2\wedge\sigma^3_1\nonumber\\
&& \qquad \qquad-\frac{7}{3}I_1\xi^2\wedge\sigma^3_1 ,\nonumber
\end{eqnarray}
where
\begin{eqnarray}\label{I}
 &&I_1=\frac{w}{\sqrt[3]{Cv^2}} ,\nonumber\\
   && I_2=-\frac{Bw+Cuw+v_t}{Cv^2} ,\\
    &&I_3=\frac{A}{Cv} .\nonumber
\end{eqnarray}
are invariants of the symmetry group of an equation of
$\mathcal{S}_2$.\\
Consider the subclass of equations (\ref{P}) such that
$Q_{v^2}\neq 0$ and $Q_{uv}\neq 0$.We denote this subclass by
$\mathcal{S}_3$. For an equation from $\mathcal{S}_3$ we normalize
$a^3_3=\frac{(Q_{v^2})^4}{(Q_{uv})^3}$ and
$b^2_2=\frac{Q_{uv}}{Q_{v^2}}$. After the absorption of torsion we
have the coframe
$\theta=\{\theta_1,\theta_2,\theta_3,\xi^1,\xi^2,\sigma^1_1,\sigma^2_1,\sigma^3_1\}$,
with the structure equations
\begin{eqnarray}\label{ES3}
&&d\theta^1=(L_1-2L_2)\theta^1\wedge\theta^2+(L_3+L_4-2I_5)\theta^1\wedge\xi^1-\theta^2\wedge\xi^2\nonumber\\
&&\qquad \quad+(L_6-2L_7)\theta^1\wedge\xi^2+\xi^1\wedge\sigma^1_1,\nonumber\\
&&d\theta^2=(3L_1-2L_8)\theta^1\wedge\theta^2+(2L_4+L_3-3L_5)\theta^2\wedge\xi^1-\theta_3\wedge\xi^2\nonumber\\
&&\qquad \quad+(2L_6-3L_7)\theta^2\wedge\xi^2+\xi^1\wedge\sigma^2_1,\nonumber\\
&&d\theta^3=(4L_1-3L_8)\theta^1\wedge\theta^3+(4L_2-3L_1)\theta^2\wedge\theta^3+(3L_6-4L_7)\theta^3\wedge\xi^2\nonumber\\
&& \qquad \quad +(3L_4+L_3-4L_5)\theta^3\wedge\xi^1+\xi^1\wedge\sigma^3_1+\xi^2\wedge\sigma^1_1,\nonumber\\
&&d\xi^1=3(L_8-L_1)\theta^1\wedge\xi^1+3(L_1-L_2)\theta^2\wedge\xi^1+3(L_7-L_6)\xi^1\wedge\xi^2,\nonumber\\
&&d\xi^2=\xi^1\wedge\theta^1+(L_8-L_1)\theta^1\wedge\xi^2-\theta^2\wedge\xi^1+(L_1-L_2)\theta^2\wedge\xi^2\nonumber\\
&&\qquad \quad+(L_4-L_5+L_9)\xi^1\wedge\xi^2,\\
&&d\sigma^1_1=\pi_3\wedge\xi^1-L_{10}\theta^1\wedge\xi^2+(5L_1-4L_8)\theta^1\wedge\sigma^1_1-L_9\theta^2\wedge\xi^2\nonumber\\
&&\qquad \quad+(5L_2-4L_1)\theta^2\wedge\sigma^1_1+(5L_5-L_3-4L_4-L_9)\xi^1\wedge\sigma^1_1\nonumber\\
&&\qquad \quad+\xi^2\wedge\sigma^2_1+(5L_7-4L_6)\xi^2\wedge\sigma^1_1,\nonumber\\
&&d\sigma^2_1=\pi_1\wedge\xi^1+L_{11}\theta^1\wedge\theta^2+\theta^1\wedge\theta^3+(6L_1-5L_8)\theta^1\wedge\sigma^2_1+\theta^2\wedge\theta^3\nonumber\\
&&\qquad \quad-L_{10}\theta^2\wedge\xi^2+(6L_2-5L_1)\theta^2\wedge\sigma^2_1-L_9\theta^3\wedge\xi^2\nonumber\\
&&\qquad \quad+(6L_7-5L_6)\xi^2\wedge\sigma^2_1+\xi^2\wedge\sigma^3_1,\nonumber\\
&&d\sigma^3_1=\pi_2\wedge\xi^1+\pi_3\wedge\xi^2+L_{11}\theta^1\wedge\theta^3+(7L_1-6I_8)\theta^1\wedge\sigma^3_1+\theta^1\wedge\sigma^1_1\nonumber\\
&&\qquad \quad +\theta^2\wedge\theta^3+(7L_2-6L_1)\theta^2\wedge\sigma^3_1+\theta^2\wedge\sigma^1_1-L_{10}\theta^3\wedge\xi^2\nonumber\\
 &&\qquad \quad+(7L_7-6L_6)\xi^2\wedge\sigma^3_1.\nonumber
\end{eqnarray}
Where the following functions
\begin{eqnarray}\label{L}
 &&L_1=\frac{Q_{uv^2}Q_{uv}}{(Q_{v^2})^3},\qquad L_7=\frac{(wQ_{v^3}+vQ_{uv^2})}{Q_{uv}},\nonumber\\
 &&L_2=\frac{Q_{v^3}(Q_{uv})^2}{(Q_v^2)^4},\hspace{.5cm} L_{10}=\frac{(Q_{v^2})^4(wQ_{uv}+vQ_{u^2})}{(Q_{uv})^4} ,\nonumber\\
 &&L_3=\frac{Q_u(Q_v)^2}{(Q_{uv})^3},\qquad L_6=\frac{Q_{v^2}(uQ_{u^2v}+vQ_{uv^2})}{(Q_{uv})^2},\\
 &&L_8=\frac{Q_{u^2v}}{(Q_{v^2})^2},\hspace{1cm} L_9=\frac{(Q_{v^2})^3(wQ_{v^2}+vQ_{uv})}{(Q_{uv})^3},\nonumber\\
 &&L_{11}=\frac{Q_{v^2}Q_{u^2}}{(Q_{uv})^2},\quad L_4=\frac{Q_{v^2}(uQ_vQ_{u^2v}+u_tQ_{u^2v}+wQ_vQ_{uv^2}+v_tQ_{uv^2})}{(Q_{uv})^4},\nonumber\\
 &&L_5=\frac{(Q_{v^2})^2(vQ_vQ_{uv^2}+u_tQ_{uv^2}+wQ_vQ_{v^3}+v_tQ_{v^3})}{(Q_{uv})^3} .\nonumber
 \end{eqnarray}
are invariants of the symmetry group of an equations of
$\mathcal{S}_3$.\\
Finally, we denote by $\mathcal{S}_4$ the subclass of equations
(\ref{P}) such that $Q_{u^2}\neq 0$ , $Q_{uv}\neq 0$ and $Q_{v^2}=
0$, modified and generalized KDV equations belong to this
subclass. For an equation from $\mathcal{S}_4$ we normalize
$a^3_3=\frac{(Q_u)^5}{(vQ_{v^2})^4}$ and
$b^2_2=\frac{uQ_{u^2}}{Q_u}$. After absorption of torsion, we have
the coframe
$\theta=\{\theta_1,\theta_2,\theta_3,\xi^1,\xi^2,\sigma^1_1,\sigma^2_1,\sigma^3_1\}$,
with the following structure equations
\begin{eqnarray}\label{ES4}
&&d\theta^1=2M_1\theta_1\wedge\theta_2+(M_4-3M_2+2M_3)\theta_1\wedge\xi_1+(2M_5-3M_6)\theta_1\wedge\xi_2\nonumber\\
&&\qquad \quad-\theta_2\wedge\xi_2+\xi_1\wedge\sigma^1_1,\nonumber\\
&&d\theta^2=(4M_1-3M_7)\theta_1\wedge\theta_2+(M_4-4M_2+3I_3)\theta_2\wedge\xi_1-\theta_3\wedge\xi_2\nonumber\\
&&\qquad \quad+(3M_5-4M_6)\theta_2\wedge\xi_2+\xi_1\wedge\sigma^2_1,\nonumber\\
&&d\theta^3=(5M_1-4M_7)\theta_1\wedge\theta_3-4M_1\theta_2\wedge\theta_3+(M_4-5M_2+4M_3)\theta_3\wedge\xi_1\nonumber\\
&& \qquad \quad -\xi_1\wedge\sigma^3_1+(4M_5-5M_6)\theta_3\wedge\xi_2+\xi_2\wedge\sigma^1_1,\nonumber\\
&&d\xi^1=3(M_7-M_1)\theta_1\wedge\xi_1+3M_1\theta_2\wedge\xi_1+3(M_6-M_5)\xi_1\wedge\xi_2,\nonumber\\
&&d\xi^2=\xi_1\wedge\theta_1+(M_7-M_1)\theta_1\wedge\xi_2+(M_8-M_2+M_3)\xi_1\wedge\xi_2,\\
&&\qquad \quad+M_1\theta_2\wedge\xi_2,\nonumber\\
&&d\sigma^1_1=\pi_3\wedge\xi_1-(M_8+M_9)\theta_1\wedge\xi_2+(6M_1-5M_7)\theta_1\wedge\sigma^1_1+\xi_2\wedge\sigma^2_1\nonumber\\
&&\qquad \quad-5M_1\theta_2\wedge\sigma^1_1+(6M_2-5M_3-M_4-M_8)\xi_1\wedge\sigma^1_1\nonumber\\
&&\qquad \quad-M_8\theta_2\wedge\xi_2-(6M_6-5M_5)\xi_2\wedge\sigma^1_1,\nonumber\\
&&d\sigma^2_1=\pi_1\wedge\xi_1-\theta_1\wedge\theta_2+\theta_1\wedge\theta_3+(7M_1-6M_7)\theta_1\wedge\sigma^2_1+\xi_2\wedge\sigma^3_1\nonumber\\
&&\qquad \quad-(M_8+M_9)\theta_2\wedge\xi_2-6M_1\theta_2\wedge\sigma^2_1)-M_8\theta_3\wedge\xi_2\nonumber\\
&&\qquad \quad+(7M_6-6M_5)\xi_2\wedge\sigma^2_1,\nonumber\\
&&d\sigma^3_1=\pi_2\wedge\xi_1+\pi_3\wedge\xi_2+\theta_1\wedge\theta_3+\theta_1\wedge\sigma^1_1+(8M_1-7M_7)\theta_1\wedge\sigma^3_1\nonumber\\
&&\qquad \quad +\theta_2\wedge\theta_3-7M_1\theta_2 \wedge\sigma^3_1-(M_8+M_9)\xi^2\wedge\sigma^3_1\nonumber\\
 &&\qquad \quad+(8M_6-7M_5)\xi^2\wedge\sigma^1_1.\nonumber
\end{eqnarray}
Where the functions
\begin{eqnarray}\label{M}
&&M_1=\frac{Q_{u^2v}(Q_{u^2})^2}{(Q_{uv})^4}, \qquad M_2=\frac{Q_{u^2v}(Q_{uv})^2(vQ_v+u_t)}{(Q_{u^2})^3},\nonumber\\
 &&M_3=\frac{(Q_{uv})^3(v_tQ_{u^2v}+vQ_vQ_{u^3}+wQ_vQ_{u^2v}+u_tQ_{u^3})}{(Q_{u^2})^4},\nonumber\\
&&M_4=\frac{Q_u(Q_{uv})^3}{(Q_{u^2})^3}, \qquad M_5=\frac{Q_{uv}(vQ_{u^3}+wQ_{u^2v})}{(Q_{u^2})^2}, \\
 &&M_6=\frac{vQ_{u^2v}}{Q_{u^2}},\qquad M_7=\frac{Q_{u^2}Q_{u^3}}{(Q_{uv})^3}, \qquad M_8=\frac{v(Q_{uv})^4}{(Q_{u^2})^3}, \nonumber\\
&&M_9=\frac{w(Q_{uv})^5}{(Q_{u^2})^4},\nonumber
\end{eqnarray}
are invariants of the symmetry group of an equation from
$\mathcal{S}_4$.\\
 The structure
equations (\ref{ES2}),(\ref{ES3}) and(\ref{ES4}) do not contain
any torsion coefficient depending on the group parameters.Their
degree of indeterminacy $r^{(1)}$ is 3, whereas the reduced
characters are $s_1=3,s_2= \ldots = s_8 = 0$. So, Cartan's test
for each of them is satisfied and the coframes are involutive.By
applying Theorem 15.12 of \cite{EIS} to above calculations we have
following statement:
\begin{thm}
The class of equation (\ref{a}) is divided into four subclasses
$\mathcal{S}_1$ to $\mathcal{S}_4$ invariant under an action of
the pseudo-group of contact transformations:\\
$\mathcal{S}_1$ consists of all systems (\ref{a}) such that $Q_{u^2}= 0$ , $Q_{u_x^2}=0$ and $Q_{uu_x}= 0$;\\
$\mathcal{S}_2$ consists of all systems (\ref{a}) such that $Q_{u^2}=0, Q_{u_x^2}= 0$ and $Q_{uu_x}\neq 0$;\\
$\mathcal{S}_3$ consists of all systems (\ref{a}) such that $Q_{u_x^2}\neq 0$ and $Q_{uu_x}\neq 0$;\\
$\mathcal{S}_4$ consists of all systems (\ref{a}) such that $Q_{u^2}\neq 0$ , $Q_{uu_x}\neq 0$ and $Q_{u_x^2}= 0$.\\
 All equation from $\mathcal{S}_1$ is equivalent to $u_{xxx}=u_t$.\\
The basic differential invariants for equations from the subclass
$\mathcal{S}_2$ are the functions $I_1, I_2$ and $I_3$ defined by
(\ref{I}). Two equations from $\mathcal{S}_2$ are equivalent with
regard to the pseudo-group of contact transformations whenever
they have the same functional dependence among the invariants
$I_1, I_2$ and
$I_3$.\\
The basic differential invariants for equations from the subclass
$\mathcal{S}_3$ are the functions $L_1,\ldots ,L_{11}$ defined by
(\ref{L}). Two systems from $\mathcal{S}_3$ are equivalent with
respect to the pseudo-group of contact transformations whenever
they have the same functional dependence among the invariants
$L_1,\ldots ,L_{11}$.\\
The basic differential invariants for equations from the subclass
$\mathcal{S}_4$ are the functions $M_1,\ldots ,M_9$ defined by
(\ref{M}). Two systems from $\mathcal{S}_4$ are equivalent with
respect to the pseudo-group of contact transformations whenever
they have the same functional dependence among the invariants
$M_1,\ldots ,M_9$.
\end{thm}

{\bf Conclusion.} In this paper, the moving coframe method of
\cite{FO1,MC} is applied to the local equivalence problem for a
class of systems of KDV-type equations under The action of a
pseudo-group of contact transformations. We found four  subclasses
and showed that every system of KDV-type equations can be
transformed to a system from one of these subclasses. The
structure equations and the invariants for all subclasses were
also found.

--------------------------------------------------------------

\end{document}